\documentclass[12pt]{amsart}
\usepackage[english]{babel}
\usepackage{amsmath}
\usepackage{amsthm}
\usepackage{amssymb}
\usepackage{mathrsfs}
\usepackage{mathdots}
\usepackage{enumerate,MnSymbol}
\usepackage[notcite, final, notref]{showkeys}
\usepackage{dsfont}
\usepackage{color}
\usepackage{shadow}
\newtheorem{theorem}{Theorem}[section]

\newtheorem{lemma}[theorem]{Lemma}
\newtheorem{proposition}[theorem]{Proposition}
\newtheorem{corollary}[theorem]{Corollary}
\newtheorem{definition}[theorem]{Definition}

\theoremstyle{definition}
\newtheorem{remark}[theorem]{Remark}

\newcommand{\BI}{\mathbb{I}}

\newcommand{\BN}{\mathbb{N}}

\newcommand{\cF}{\mathcal{F}}

\newcommand{\cH}{\mathcal{H}}

\newcommand{\g}{\mathbf{g}}

\newcommand{\ov}{\overline}

\newcommand{\sesq}[1]{\left[  #1 \right] }

\newcommand{\qexp}[1]{e_q^{#1}\, }
\newcommand{\qln}[1]{\ln_q #1\, }

\usepackage{hyperref}

\title[Tsallis' $q$-analysis]{Tsallis' $q$-analysis, new scales of interpolating spaces, and $q$-rational functions}

\author[D. Alpay]{Daniel Alpay}
\address{(DA)
Faculty of Mathematics, Physics, and Computation\\
Schmid College of Science and Technology\\
Chapman University\\
One University Drive
Orange, California 92866\\
USA}
\email{alpay@chapman.edu}

\author[P. Cerejeiras]{Paula Cerejeiras}
\address{(PC) Center for research and development in mathematics and applications\\Department of Mathematics, University of Aveiro  \\
Campus Universit\'ario de Santiago  \\ 3810-193 Aveiro\\Portugal}
\email{pceres@ua.pt}

\author[U. Kaehler]{Uwe Kaehler}
\address{(UK) Center for research and development in mathematics and applications\\Department of Mathematics, University of Aveiro  \\
Campus Universit\'ario de Santiago  \\ 3810-193 Aveiro\\Portugal}
\email{ukaehler@ua.pt}

\begin{document}
\maketitle
\begin{abstract}
In this paper we are studying the fundamental tools for a quantum calculus based on the Tsallis $q$-exponential In particular we are looking at $q$-Fock spaces, structural identities, as well as rational functions in this context. 
\end{abstract}

\tableofcontents

\section{Introduction}

The exponential function $e^x$ plays a key role in Hida's white noise space theory (see e.g. \cite{Hida_BM,MR1244577,MR2444857,MR1408433})
because of two complementary properties: first the coefficients of
its power series expansion at the origin are positive, and next the function
\[
  e^{-\frac{\|f-\g\|^2}{2}}
\]
is positive definite, in the sense of reproducing kernels, for $f,g$ in the space of real Schwartz functions $\mathcal S_{\mathbb R}$.
The second property allows to  construct, thanks to the Bochner-Minlos theorem, the underlying probability space $(\mathcal S_{\mathbb R}^\prime,\mathcal C,P)$,
where $\mathcal C$ is the sigma-algebra generated by the cylinder sets. The first property allows to show
the unitary equivalence between $\mathbf L_2(\mathcal S_{\mathbb R}^\prime,\mathcal C,P)$ and the reproducing kernel Hilbert space with reproducing
kernel $e^{\sum_{n=1}^\infty z_n\overline{w_n}}$.

These properties and the example of the Mittag-Leffler
function lead us in \cite{stochastics_aveiro_1} to introduce a class of functions, denoted by $ML$, for which Hida's theory can be extended. The Mittag-Leffler function
belongs to $ML$ and corresponds to the grey noise theory developed by Jahnert \cite{Jahnert}.\smallskip

In a completely different avenue, the exponential function plays a key role in equilibrium statistical physics via the Boltzmann-Gibbs distribution.
It also plays a key role to define the entropy via its inverse version, namely $\ln x$.\smallskip

The exponential function is the unique solution of the ODE $y^\prime=y$ with initial condition $y(0)=1$.
Tsallis extended the domain of equilibrium statistical physics to a domain called {\sl Nonextensive Statistical Mechanics}, replacing this ODE by
\begin{equation}
y^\prime =y^q,\quad y(0)=1,
  \end{equation}
  where $q\in[1,\infty)$ is pre-assigned; see e.g.
  \cite[(11) p. 458]{Tsallis-triplet}. The solution
\[
e_q(z)=(1+(1-q)x)^{\frac{1}{1-q}}
  \]
  now plays the role of the exponential, with the logarithm being

  \[
\ln_q x=\frac{x^{1-q}-1}{1-q}.
    \]
   At this stage we remark that $\ln_q$ is not the usual quantum logarithm. Note that
    \[
x^q\ln_qx=x\ln_{2-q}x=\frac{x-x^q}{1-q}.
      \]
      For a probability distribution $p_1,\ldots, p_N$ the $S_q$ entropy is defined by
      \[
S_q=k\frac{1-\sum_{n}^Np_n^q}{q-1}=-k\left(\sum_{n}^Np_n^q\ln_qp_n\right)
        \]
        corresponding to the classical entropy for $q=1$. Note the difference with the Renyi entropy
        \[
          S_{\rm Renyi}=\frac{\ln(\sum_{n=1}^Np_n^q)}{1-q};
        \]
        See e.g. \cite[p. 539]{Tsallis-triplet} for the latter.

  To make links with reproducing kernel spaces, we already remark that in the limiting case $q=1$ the function
  $e_1(z\overline{w})=e^{z\overline{w}}$ (corresponding to the Fock space) while for $q=2$, we have $e_2(z\overline{w})=\frac{1}{1-z\overline{w}}$,
  corresponding to the Hardy space. So from the analytic point of view, the parameter $q$ allows a natural interpolation between the Hardy space and the Fock space.

  Such parametrizations already appear for the Mittag-Leffler function, with associated reproducing kernel $E_\alpha(n\alpha+1)$, with $\alpha=1$
  corresponding to the Fock space and $\alpha=0$ to the Hardy space, and in $q$-analysis; see \cite{Cerejeiras23} for the latter.\\

  The parameter $q$ relates to a wide range of different settings, such as the entropy of hydrogen atoms, wavelet analysis
and the volatility of financial markets; see \cite{nelson_kendric}.

\section{Preliminaries} 
The Fock space  is the reproducing kernel Hilbert space with reproducing kernel
  \[
e^{z\overline{w}}=\sum_{n=0}^\infty \frac{z^n\overline{w}^n}{n!},\quad z,w\in\mathbb C,
    \]
    and plays a key role in quantum mechanics, while the Hardy space of the open unit disk $\mathbb D$ is the reproducing kernel Hilbert space with reproducing kernel
\[
  \frac{1}{1-z\overline{w}}=\sum_{n=0}^\infty z^n\overline{w}^n, \quad z,w\in\mathbb D,
\]
and plays a key role in operator theory and linear system theory.  \\

The Hardy space can be characterized (up to a positive multiplicative factor for the inner product) as the Hilbert space of power series converging at the origin and such that
\begin{equation*}
  R_0^*=M_z,
\end{equation*}
  where $M_z$ is the operator of multiplication by $z$ and
\begin{equation*}
  R_0f(z)=\frac{f(z)-f(0)}{z}.
  \end{equation*}

  Both the Fock and the Hardy space are reproducing kernel Hilbert spaces of power series and have been extended and generalized in numerous directions, of which we mention in particular the $q$-quantum calculus approach; see \cite{Cerejeiras23}. In the present work we focus on another $q$ approach based on Tsallis work on statistics; see \cite{Yamano2002}.  \smallskip

Following \cite{Yamano2002}, who considers the case $q \in (0,1),$ we introduce the $q$-logarithm and $q$-exponential defined by Tsallis \cite{Tsallis2009},
\begin{equation} \label{Tsallis_log}
\qln{x} := \frac{x^{1-q}-1}{1-q}, \quad x >0,
\end{equation}
and
\begin{equation} \label{Tsallis_exp}
\qexp{x} := \left\{  \begin{array}{cl}
[1+(1-q)x]^{\frac{1}{1-q}}, & 1+(1-q)x \geq 0 \\
 & \\
 0, & {\rm otherwise.}
\end{array}  \right. .
\end{equation}

We observe that these functions are the inverse of each other. When $q \to 1^-$ then we have that they coincide with the classical logarithm and exponential functions, i.e.
$$\lim_{q \to 1^-} \qln{x} = \ln x, \qquad \lim_{q \to 1^-} \qexp{x} = e^x.$$
When $q \to 0^+$ we get the degenerated case
$$\lim_{q \to 0^+} \qln{x} = x-1, \qquad \lim_{q \to 0^+} \qexp{x} = x+1.$$

As in the case of the Fock space this exponential generates a reproducing kernel $K_q(z,w) := \qexp{z\ov w} = [1+(1-q)z\ov w ]^{\frac{1}{1-q}}$ where $q\not=1.$ In addition, if $\frac{1}{1-q} = k \in \BN,$ i.e. $q = \frac{k-1}{k},$ we get a polynomial kernel of degree $k,$
\begin{equation} \label{RK}
K_{\frac{k-1}{k}}(z,w) = \Big( 1+\frac{z\ov w}{k}  \Big)^k, \quad k=1, 2, 3, \ldots,
\end{equation} which is a positive definite kernel giving rise to a finite-dimensional reproducing kernel Hilbert space. For example, for $q= \frac{3}{4}$ we have the 4th degree polynomial
$$K_{\frac{3}{4}}(z,w) =\Big( 1+\frac{z\ov w}{4}  \Big)^4= 1+ z\ov w + \frac{3}{8} (z\ov w)^2 +  \frac{1}{16} (z\ov w)^3 + \frac{1}{256} (z\ov w)^{4}.$$

For the remaining cases of $q \in (0,1)$ we obtain the reproducing kernel in terms of a binomial power series expansion
\begin{eqnarray}
K_q(z,w) & = & \qexp{z\ov w} \nonumber \\
 & = & \big[ 1+(1-q)z\ov w \big]^{\frac{1}{1-q}} \nonumber \\
&:=& \sum_{k=0}^\infty \frac{(1/ (1-q))_k}{k !} [(1-q) z\ov w]^k, \qquad |z\ov w| < \frac{1}{1-q}, ~q \in (0,1)  \nonumber \\
&=& \sum_{k=0}^\infty \Big[(1-q)^k \,\Big(\frac{1}{1-q}\Big)_k \Big] \frac{(z\ov w)^k}{k !}, \label{qRK} 
\end{eqnarray}
whereas $(a)_0=1, ~(a)_k = a (a-1) \cdots (a-k+1),$ denotes the lowering factorial. Easy computations show that $(1-q)^0 \,\Big(\frac{1}{1-q}\Big)_0 = (1-q)^1 \, \Big(\frac{1}{1-q}\Big)_1 =1,$ while for $k\geq 2,$
$$
\alpha_k = (1-q)^k \,\Big(\frac{1}{1-q}\Big)_k = (1-(1-q)) (1-2 (1-q)) \cdots \big(1-(k-1)(1-q) \big).
$$
In terms of $q$ this corresponds to the following coefficients in (\ref{qRK}):
\begin{equation}\label{coeff}
\alpha_0=\alpha_1=1, \quad \alpha_k = q(2q-1)(3q-2) \cdots \big( (k-1)q -(k-2) \big), \quad k=2, 3, \ldots.
\end{equation}
Solving the diophantine equation $(k-1)q -(k-2) =0$ we conclude that  $$\alpha_k = (1-q)^k \,\Big(\frac{1}{1-q}\Big)_k >0$$ for $k=0,1, \ldots, k_0,$ while $\alpha_k = (1-q)^k \,\Big(\frac{1}{1-q}\Big)_k$ alternates sign for $k>k_0,$ with $k_0 = \lfloor \frac{2-q}{1-q} \rfloor \geq 2.$ This will imply the appearance of Krein spaces in the sequel for $q\in(0,1)$.

 Furthermore, we can consider the above series  (\ref{RK}) also in the case $q\in (1,2]$.

\section{Reproducing kernel spaces} A careful look into (\ref{Tsallis_log}) and (\ref{Tsallis_exp}) shows that the parameter $q \in (0,1)$ can be extended to $q >0,$ pending a correct adaptation for $q=1$ and imposing $|z\ov w| < \frac{1}{|1-q|}$ in (\ref{qRK}). Furthermore, we remark that $q=2$ corresponds to
$$\lim_{q \to 2} \qln{x} = 1-\frac{1}{ x}, \qquad \lim_{q \to 2} \qexp{x} = \frac{1}{1-x},$$
for $0<x <1.$ This last case is linked to the reproducing kernel of the Hardy space.\\

For each $q>0$ we have the following sesquilinear form (or indefinite inner product) defined on the monomials as
\begin{equation} \label{Inner_Fock}
\sesq{z^n, z^k}  = \frac{k!}{\alpha_k} \delta_{n,k} =  \left\{  \begin{array}{cl}
\delta_{n,k}, & k=0,1 \\
~& \\
\frac{k!}{q(2q-1)(3q-2) \cdots \big( (k-1)q -(k-2) \big)} \delta_{n,k}, & k \geq 2
\end{array}
\right. ,
\end{equation} and extended by linearity to the whole space under the condition that $nq \not= n-1,$ for all $n \in \BN.$ In this case we obtain a generalized Fock space with reproducing kernel
\begin{equation} \label{RKHS}
K_q(z, w) =  \qexp{z\overline w}  = 1+z\overline w +\sum_{k=2}^\infty \frac{\Pi_{n=1}^{k-1} [nq -(n-1)] }{k!} (z\overline w)^k, \quad |z\ov w| < \frac{1}{|q-1|}.
\end{equation} When $q=1$ we obtain $K_1(z, w) =  e^{z\overline w},$ for $|z\ov w| < +\infty,$ the reproducing kernel for the classical Fock space. For $q>1$ we have $\Pi_{n=1}^{k-1} [nq -(n-1)] >0$ for all $k$. For $q \in (0,1)$ the product $\Pi_{n=1}^{k-1} [nq -(n-1)]$ has alternating sign for  $k > k_0 := \lfloor \frac{1}{1-q} \rfloor,$ while it is non-negative for $k \leq k_0.$
This leads to the following lemma.

\begin{lemma} \label{Lem:1} For $q \in (0,1)$ such that  $nq \not= n-1,$ for all $n \in \BN$, the reproducing kernel (\ref{RKHS}) corresponds to a Krein space. However, when $q = \frac{n-1}{n}$ for some $n \in \BN,$ then the reproducing kernel (\ref{RKHS}) has only a finite number of positive terms. In this case, the space is a finite dimensional Hilbert space.
\end{lemma}


\subsection{Krein spaces}
We first give the general definition of a Krein space, and then study Krein spaces of power series. We consider the case of complex numbers as coefficients, but the case of real numbers or of quaternions are also of interest.
Given an Hermitian form on a complex vector space $\mathcal V$ we take the convention
\begin{equation}
[\lambda u,\mu v]=\overline{\mu}\lambda[u,v],\quad u,v\in\mathcal V,\quad \lambda,\mu\in\mathbb C
  \end{equation}
for the linearity and anti-linearity of the coefficients.

\begin{definition}
  A complex vector space $\mathcal V$  space endowed with an Hermitian form $[\cdot,\cdot]$ is called a Krein space if it can written as a direct and orthogonal sum as
  \begin{equation}
    \label{funda-decomp}
\mathcal V=\mathcal V_+[+]\mathcal V_-,
\end{equation}
where $(\mathcal V_+,[\cdot,\cdot])$ and $(\mathcal V_-,-[\cdot,\cdot])$ are Hilbert spaces, such that
\[
  \mathcal V_+\cap\mathcal V_-=\left\{0\right\}\quad (direct \,\,sum)
\]
and
\[
  [v_+,v_-]=0,\quad \forall v_+\in\mathcal V_+\,\, and\,\,\forall v_-\in\mathcal V_-\quad (orthogonal\,\,sum).
\]
  \end{definition}

  \eqref{funda-decomp} is called a fundamental decomposition. It is not unique, but in  degenerate cases, when $(\mathcal V,[\cdot,\cdot])$ is either a
  Hilbert space or an anti-Hilbert space.\smallskip

  Given a fundamental decomposition \eqref{funda-decomp}, consider the Hermitian form
  \[
    \langle v,w\rangle=[v_+,w_+]-[v_-,w_-]
  \]
  where $v=v_++v_-$ and $w=w_++w_-$ are the corresponding decompositions of $v,w\in\mathcal V$. Then $(\mathcal V,\langle\cdot,\cdot\rangle)$ is
  a Hilbert space with associated norm $\sqrt{[v_+,v_+]-[v_-,v_-]}$. The norm depends on the decomposition (but in degenerate cases, when it is unique),
  but all the norms obtained are equivalent; see \cite[p. 102]{bognar}.\smallskip

  Let $A\subset \mathbb N_0$ and let $(\gamma_a)_{a\in A}$ be a family of real numbers different from $0$.
Let
\[
  A_+=\left\{a\in A\,\,;\,\, \gamma_a>0\right\}\quad{\rm and}\quad   A_-=\left\{a\in A\,\,;\,\, \gamma_a<0\right\}.
\]
Suppose furthermore that the power series $\sum_{a\in A}\frac{z^n}{\gamma_a}$ has a strictly positive radius of convergence, i.e.
\[
R=\frac{1}{\limsup |1/ \gamma_a|^{1/a}}=\liminf |\gamma_a|^{1/a}>0.
  \]

  \begin{proposition}
    In the above notation and with $R>0$, the space of functions
    \begin{equation}
\mathfrak K(\gamma)=\left\{f : f(z)=\sum_{a\in A}f_az^a \mbox{~is convergent for }|z| < R \right\},\end{equation}
with Hermitian sesquilinear form
\begin{equation}
[f,g]=\sum_{a\in A}=\gamma_a\overline{g_a}f_a,\qquad (g(z)=\sum_{a\in A}g_az^a\in\mathfrak K),
\end{equation}
is a reproducing kernel Krein space with reproducing kernel
\begin{equation} \label{eq:3.12}
  K(z,w)=\sum_{a\in A}\frac{z^a\overline{w}^a}{\gamma_a}.
  \end{equation}
\end{proposition}

\begin{proof}
  In the above notation, define
  \[
    [ f, g ]=\sum_{a\in A} \gamma_a \,\overline{g_a}f_a
  \]
  and
  \[
    (Jf)(z)=\sum_{a\in A_+}f_az^a-\sum_{a\in A_-}f_az^a.
  \]
  Then, $(\mathfrak K(\gamma),\langle\cdot,\cdot\rangle)$ is a Hilbert space and $J$  is a bounded map from this Hilbert space into itself, satisfying $J=J^*=J^{1}$.
  Furthermore
  \begin{equation}
    [f,g]=\langle f, Jg\rangle
    \end{equation}
      and this ends the proof.
  \end{proof}


\begin{definition} The  $q$-Fock-Tsallis space $\cF_q$ is given by
\begin{equation}\label{F_q}
\cF_q = \{ f = \sum_{k=0}^\infty f_k z^k: \quad  \| f\|_{\cF_q}^2 =  \sum_{k=0}^\infty  |f_k|^2 |\gamma_k| < +\infty  \}.
\end{equation}
\end{definition}

\begin{theorem} The $q$-Fock-Tsallis space $\cF_q$ has reproducing kernel $K_q$ associated to the sequence of (non-zero) real numbers $(\gamma_k)_{k=0}^\infty,$ where $\gamma_k= \frac{k!}{ \alpha_k} ,$ with $\alpha_k = (1-q)^k \,\big(\frac{1}{1-q}\big)_k$, i.e.
\begin{equation}\label{H(K_q)}
 f(z) = \sesq{f, K_q(\cdot,z)}, \quad \mbox{for all }f \in \cF_q.
\end{equation}
 \end{theorem}

  \begin{proof} Following (\ref{eq:3.12}) and using the given family of (non-zero) real numbers $$\gamma_k=\frac{k!}{(1-q)^k \,\big(\frac{1}{1-q}\big)_k}, \quad k=0,1, \cdots,$$ we have the $q$-Fock-Tsallis space defined as
   $$\cF_q = \{ f = \sum_{k=0}^\infty f_k z^k: \quad  \| f\|_{\cF_q}^2 =  \sum_{k=0}^\infty  |f_k|^2 \frac{k!}{ |1-q|^k \,\big|\big(\frac{1}{1-q}\big)_k \big|}  < +\infty  \}.$$
   Since $K_q(w,z) = \qexp{w\overline z}$ (see (\ref{RKHS})) we get for $f = \sum_{k=0}^\infty f_k z^k \in \cF_q$
   \begin{gather*}
   \sesq{f, K_q(\cdot,z)} = \sum_{k=0}^\infty \frac{k!}{ (1-q)^k \,\big(\frac{1}{1-q}\big)_k} \, \overline{(1-q)^k \,\big(\frac{1}{1-q}\big)_k \frac{\ov z^k}{k!}} \, f_k  = \sum_{k=0}^\infty  f_k z^k = f(z).
   \end{gather*}

    \end{proof}

De facto, the $q$-Fock-Tsallis space is a reproducing kernel Hilbert space whenever $q \geq 1$ since in this case we have $\alpha_k = (1-q)^k \,\Big(\frac{1}{1-q}\Big)_k >0,$ for all $k$. In the case $q\in (0,1)$  It is still a reproducing kernel Hilbert space when we have $q= \frac{k-1}{k}$ for some $k \in \BN.$ Otherwise we have a reproducing kernel Krein space. We also remark that there is no uniqueness with respect to the Krein space, i.e. for the same reproducing kernel $K_q$ there might exist different Krein spaces associated to
    this reproducing kernel. See the paper of L. Schwartz \cite[Sections 12 and 13]{schwartz} and for a later counterexample,  \cite{a2}.


\subsection{Main operators and their properties} In such a generalized Fock space (associated to $q >0$ s.t. $nq \not= n-1,$ for all $n \in \BN$), the first operator to consider is the classical multiplication operator
\begin{equation} \label{multiplicationOp}
M_z z^n = z^{n+1}, \quad n=0,1, 2, \ldots
\end{equation}
Its adjoint $M_z^\ast$ acting on $z^k$ is given as $M_z^\ast z^k = \beta_{k-1}z^{k-1},~k=1,2, \cdots,$ which leads to
\begin{eqnarray*}
\sesq{M_z z^n, z^k} = \sesq{z^n, M_z^\ast z^k} &\Leftrightarrow& \frac{k!}{(1-q)^k \,\Big(\frac{1}{1-q}\Big)_k} \delta_{n+1,k} = \beta_{k-1} \frac{(k-1)!}{(1-q)^{k-1} \,\Big(\frac{1}{1-q}\Big)_{k-1}} \delta_{n,k-1}.
\end{eqnarray*}
We obtain $\beta_0 =1$ while for the remaining terms we get
$$\beta_{k-1} = \frac{k}{(k-1)q -(k-2)}, \quad k=2,3, \ldots.$$

We remark that $\beta_{k}'s$ are negative for $k > k_0 := \lfloor \frac{1}{1-q} \rfloor,$ and positive otherwise.

Due to its similarity with the integration operator $\BI z^n = \frac{z^{n+1}}{n+1}, ~n=0,1, 2, \ldots,$ we obtain for its adjoint $\BI^\ast z^n=\iota_{n-1}z^{n-1}, n=1,2, \ldots,$ with $\iota_{n-1}$ constant depending on $n$, that
$$\sesq{\BI z^n, z^k} = \sesq{z^n, \BI^\ast z^k} \Leftrightarrow \frac{k!}{(n+1)! (1-q)^k \,\Big(\frac{1}{1-q}\Big)_k} \delta_{n+1,k} = \iota_{k-1} \frac{(k-1)!}{(1-q)^{k-1} \,\Big(\frac{1}{1-q}\Big)_{k-1} } \delta_{n,k-1},$$
which gives $\iota_0 = 1$ and
$$\iota_{k-1}= \frac{1}{(k-1)q -(k-2)}, \quad k=2,3, \ldots$$

For the adjoint of $R_0,$ its action on $z^k$ given by $R_0^\ast z^k = \tau_{k+1}z^{k+1}, k=0,1,2,\ldots,$ we obtain
\begin{eqnarray*}
\sesq{R_0 z^n, z^k} = \sesq{z^n, R_0^\ast z^k} &\Leftrightarrow& \frac{k!}{(1-q)^k \,\Big(\frac{1}{1-q}\Big)_k} \delta_{n-1,k} = \tau_{k+1} \frac{(k+1)!}{(1-q)^{k+1} \,\Big(\frac{1}{1-q}\Big)_{k+1}} \delta_{n,k+1},
\end{eqnarray*}
so that $\tau_1 = 1$ and $$\tau_{k+1}= \frac{kq -(k-1)}{k+1},\quad k=1,2,\ldots$$

Those results summarize as follows:

\begin{lemma} \label{Lem:3} For $q >0$ such that  $nq \not= n-1,$ for all $n \in \BN$, we have the following adjoint operators of the backward-shift operator $R_0$, of the multiplication operator $M_z$, and of the operator of integration $\BI$:
\begin{eqnarray*} \label{M^ast}
R^\ast_0 z^k & = & \frac{kq -(k-1)}{k+1} z^{k+1}, \\
M_z^\ast z^k & = & \frac{k}{(k-1)q -(k-2)}z^{k-1},\\
\BI^\ast z^k & = & \frac{1}{(k-1)q -(k-2)}z^{k-1} = \frac{1}{k} M_z^\ast z^k,
\end{eqnarray*}
for  $k=1, 2,3,\ldots$.

\end{lemma}

Therefore, we obtain a scale of generalized Fock spaces with the following particular examples:
\begin{itemize}
\item for $q=2$ we obtain the Hardy space case with reproducing kernel
$$K_2(z,w) = e_2^{z \ov w} = \frac{1}{1-z \ov w},$$
and operators
$$M_z^\ast z^k = z^{k-1} = R_0 z^k \mbox{\quad and \quad }R_0^\ast z^k = z^{k+1} = M_z z^k;$$
\item  for $q=1$ we obtain the Fock space case with reproducing kernel
$$K_1(z,w) = e^{z \ov w},$$
and operators
$$M_z^\ast z^k = k z^{k-1} = \partial z^k \mbox{\quad and \quad} R_0^\ast z^k = \frac{z^{k+1}}{k+1} = \mathbb{I} z^k,$$ where we recall that $\mathbb I$ denotes the integration operator,
$$\mathbb I z^k = \frac{z^{k+1}}{k+1}, \quad k=0,1,2, \ldots.$$
\item for $q=0$ one has to consider the Fock space restricted to linear functions $f(z) = a_0+a_1z.$ The reproducing kernel is then
$$K_0(z,w) = e_0^{z \ov w} =1+z \ov w,$$
with operators
 $$M_z^\ast z^k =-\frac{k}{k-2} z^{k-1} \mbox{\quad and \quad} R_0^\ast z^k = -\frac{k-1}{k+1} z^{k+1},\qquad k=0,1.$$
\end{itemize}

\begin{remark} This scale of Fock spaces is different from the more classical $\textbf{q}$-quantum scale in~\cite{MR495971,Cerejeiras23}, also called $q$-Fock spaces or Arik-Coon spaces, where the Hardy space corresponds to $\textbf{q}=0$ and the classical Fock space corresponds to $\textbf{q}=1.$
\end{remark}

This also leads to
$$R_0^\ast z^k =  \frac{kq -(k-1)}{k+1} z^{k+1} = q z^{k+1} - (k-1+q) \frac{z^{k+1}}{k+1} = q M_z z^k - (k-1+q) \mathbb I z^k$$
$$= q M_z z^k - (k+1+q-2) \mathbb I z^k = [(q-1) M_z  -(q-2) \mathbb I ]z^k,$$
that is,
\begin{equation} \label{R_0^astexp}
R_0^\ast  = (q-1) M_z  -(q-2) \mathbb I,
\end{equation}
while for $q\not=1,$
\begin{equation} \label{M_z^astexp}
M_z^\ast  = \frac{1}{q-1}\big(R_0  +(q-2) \mathbb I^\ast \big).
\end{equation}

The above considerations can be summed up in the following theorem.

\begin{theorem}\label{RI-theorem}
  Let $q\ge 1$. The $q$-Fock-Tsallis space is the only reproducing kernel Hilbert space of power series which contains the polynomials and where it holds
  \[
    (q-1)M_z^*=R_0+(q-2)\mathbb I^*
    \]
  \end{theorem}

    \begin{remark}
      We note that $q=1$ leads to $R_0=\mathbb I^*$, which characterizes the classic Fock space, while $q=2$ gives $M_z^*=R_0$, which characterizes the Hardy space.
      \end{remark}

We conclude with the  boundedness of the multiplication operator in the $q$-Fock-Tsallis space.

\begin{lemma} \label{Lem:5} For $q >0$ such that  $kq \not= k-1,$ for all $k \in \BN$, we have the operator $M_z$ is bounded from $\cF_q$ into itself whenever $q < 1$.
\end{lemma}
\begin{proof}
For $f = \sum_{k=0}^\infty f_k z^k$ we have $M_z f = \sum_{k=0}^\infty f_k z^{k+1}.$ Hence,
$$\| M_z f\|_{\cF_q}^2 =  \sum_{k=0}^\infty  |f_k|^2 |\gamma_{k+1}| =  \sum_{k=0}^\infty \big( |f_k|^2 |\gamma_{k}| \big) \frac{|\gamma_{k+1}|}{|\gamma_k |}.$$
Now for $k >1$ we get
$$\frac{|\gamma_{k+1}|}{|\gamma_k|} = \frac{k+1}{k} \Big| \frac{q(2q-1)(3q-2) \cdots \big( (k-1)q -(k-2) \big)} {q(2q-1)(3q-2) \cdots \big( (k-1)q -(k-2) \big)\big( kq -(k-1) \big)} \Big|= \frac{k+1}{k |kq -(k-1)|},$$ so that we obtain $$\frac{|\gamma_{k+1}|}{|\gamma_k|} < 1$$ for large enough $k$ and
$$\| M_z f\|_{\cF_q}^2 \leq C \| f\|_{\cF_q}^2$$
for some $C >0$, i.e. the result holds.
\end{proof}

By Lemma \ref{Lem:1} we have that the reproducing kernel space $\mathcal H(K_q)$ for $q \in (0,1)$ associated to the reproducing kernel $K_q(z,w)$ (see (\ref{RKHS})) is a Krein space in case of $nq \not= n-1,$ for all $n \in \BN.$

Furthermore, we can state the following result.
\begin{theorem}
If $M_z$ is a bounded operator from $\mathcal H(K_q)$ into itself, then
\begin{equation} \label{M_z^ast}
M_z^\ast K_q(\cdot, w) = \ov w K_q(\cdot, w), \quad \mbox{\rm for all }w.
\end{equation}
\end{theorem}

\begin{proof} By (\ref{RKHS}) we have
\begin{eqnarray*}
M_z^\ast K_q(z, w) & = & M_z^\ast  \qexp{z\overline w} \\
& = & M_z^\ast  \left( 1+z\overline w +\sum_{k=2}^\infty \frac{\Pi_{n=1}^{k-1} [nq -(n-1)] }{k!} (z\overline w)^k \right) \\
& = &  0 + (M_z^\ast z)\overline w +  (M_z^\ast z^2) \overline w^2+ \sum_{k=3}^\infty \frac{\Pi_{n=1}^{k-1} [nq -(n-1)] }{k!} \overline w^k (M_z^\ast z^k)  \\
& = &   \overline w +  z \overline w^2+ \sum_{k=3}^\infty \frac{\Pi_{n=1}^{k-1} [nq -(n-1)] }{k!} \overline w^k \frac{k}{(n-1)q -(n-2)}z^{k-1} \quad \mbox{by Lemma \ref{Lem:3}}  \\
& = &   \overline w +  z \overline w^2+ \sum_{k=3}^\infty \frac{\Pi_{n=1}^{k-2} [nq -(n-1)] }{(k-1)!} \overline w^k z^{k-1}  \\
& = &   \overline w \left( 1+  z \overline w+ \sum_{k=2}^\infty \frac{\Pi_{n=1}^{k-1} [nq -(n-1)] }{k!} (z\overline w)^{k} \right) = \overline w  K_q(z, w), \quad \mbox{for }|z\ov w| < \frac{1}{|1-q|}.
\end{eqnarray*}

\end{proof}
This theorem is a particular case of for instance \cite[5.22 p. 78]{MR3526117} and \cite[Theorems 3.5 and 3.6]{MR3050315}. Multipliers are an important tool in reproducing
kernel spaces theory. For a early instance of use (leading to the linear fractional transformation in Schur analysis), see  \cite[p. 339]{dbr1}.\smallskip

Let us now consider a general version of Theorem~\ref{RI-theorem} which is valid in a broader context than the one considered in this paper.

  \begin{theorem}
    Let $a,b\in\mathbb C$ be such that $a+b\not=0,$ $a\neq 0$, $(\tilde \gamma_n)_{n=0}^\infty$ a sequence of (non-zero) real numbers and let $\mathfrak K(\tilde \gamma)$ be a reproducing kernel Krein space of functions in which the polynomials
    form a dense set, with reproducing kernel
    \[
K(z,w)=\sum_{n=0}^\infty\frac{z^n\overline{w}^n}{\tilde \gamma_n}
\]
where all the $\tilde \gamma_n$ are assumed different from $0$. Assume that $R_0$ is defined on the polynomials.  Then,
    \begin{equation}
      \label{aba123}
M_z^*=\frac{aR_0+b\mathbb I^*}{a+b}
\end{equation}
if and only if
\begin{equation}
  \tilde \gamma_n=\frac{\tilde \gamma_0}{\prod_{k=1}^n\left(1+\frac{b}{a}\left(1-\frac{1}{k}\right)\right)},\quad n=1,2,\ldots
  \end{equation}
    \end{theorem}

    \begin{proof}
      By hypothesis
      \begin{equation}
        [ z^n , z^m ]_{\mathfrak K}=\tilde \gamma_n
        \end{equation}
      Then, on the one hand,
        \[
          [ M_z^*z^n, z^{n-1} ]_{\mathfrak K} = [ z^n , z^n]_{\mathfrak K} = \tilde \gamma_n
        \]
        while, taking into account \eqref{aba123}, and since $R_0$ is defined on polynomials, we get
          \begin{eqnarray*}
            [ M_z^* z^n , z^{n-1} ]_{\mathfrak K} & = & [ \frac{aR_0+b\mathbb I^*{a+b}} z^n , z^{n-1} ]_{\mathfrak K} \\
            & = & \frac{a}{a+b} [ R_0z^n, z^{n-1} ]_{\mathfrak K}+\frac{b}{a+b} [ \mathbb I^* z^n , z^{n-1}]_{\mathfrak K} \\
          & = & \frac{a}{a+b} [ z^{n-1} , z^{n-1} ]_{\mathfrak K} + \frac{b}{a+b}  [ z^n,\mathbb I z^{n-1}]_{\mathfrak K}
          \\
      & = & \frac{a}{a+b} \tilde \gamma_{n-1}+ \frac{b}{a+b} [ z^n , \frac{z^n}{n}]_{\mathfrak K}\\
                                    & = & \frac{a}{a+b}\tilde \gamma_{n-1}+\frac{b}{a+b}\frac{\tilde \gamma_n}{n},
          \end{eqnarray*}
          and so
          \[
            \tilde \gamma_n=\frac{a}{a+b}\tilde \gamma_{n-1}+\frac{b}{a+b}\frac{\tilde \gamma_n}{n}\quad n=1,2,\ldots.
            \]
        \end{proof}

 For the particular case of $M^\ast_z$ and $\mathbb I^\ast$ given as in Lemma~\ref{Lem:3} (and Theorem~\ref{RI-theorem} for $q>1$) we have the following statement.

\begin{lemma} \label{Lem:tobeseen} For $q >0$ such that  $nq \not= n-1,$ for all $n \in \BN$, we have for the constants in (\ref{aba123}) the (non-unique) pair of values $a=1$ and $b=q-2$.
\end{lemma}

We remark that this pair $a=1$ and $b=q-2$ corresponds to the sequence
\begin{equation} \label{eq:3.25}
\gamma_0 =1, \quad \gamma_n = \frac{n!}{\prod_{k=1}^n\left((k-1) q -(k-2)\right)}.
\end{equation}

We also remark that the case $q=2$ leads to $b=0$ and $M_z^\ast = R_0,$ as in the case of the Hardy space.

          \begin{lemma}
            \label{FHab}
            With $\alpha=b/a$ we have
            \begin{equation}
              1+\sum_{n=1}^\infty z^n\prod_{k=1}^n\left(1+\frac{b}{a}\left(1-\frac{1}{k}\right)\right)=\begin{cases}(1-(1+\alpha)z)^{-\frac{1}{\alpha+1}},
                  \quad \alpha\not=-1\\
                  \, \,e^z,\quad \hspace{2.7cm}\alpha=-1.
                  \end{cases}
              \end{equation}
            \end{lemma}

            \begin{proof}

          A direct calculation gives
              \[
                \begin{split}
                  1+\sum_{n=1}^\infty z^n\prod_{k=1}^n\left(1+\frac{b}{a}\left(1-\frac{1}{k}\right)\right)&=
                  1+\sum_{n=1}^\infty z^n\prod_{k=1}^n\left(                    k+\alpha k-\alpha\right)\\
                  &=                  1+\sum_{n=1}^\infty((1+\alpha) z)^n\prod_{k=1}^n\left(                    k-\frac{\alpha}{\alpha+1}\right)\\
                  &=                  1+\sum_{n=1}^\infty(-(1+\alpha) z)^n\prod_{k=1}^n\left(                    \frac{\alpha+1-1}{\alpha+1}-k\right)\\
                                    &=                  1+\sum_{n=1}^\infty(-(1+\alpha) z)^n\prod_{k=1}^n\left(           -         \frac{1}{\alpha+1}+1-k\right).
                \end{split}
              \]
              Since
              \[
(1+z)^{-\frac{1}{\alpha+1}}=1+\sum_{n=1}^\infty z^n\prod_{k=1}^n\left(           -         \frac{1}{\alpha+1}+1-k\right).
                \]
                we have our result.
              \end{proof}
    In the case $a=-b$ where the above result does not hold we have the classic Fock space. For the remaining cases $a \not= -b$ we obtain the sequence
    $$\tilde \gamma^{a,b}_0= 1, \quad \tilde \gamma^{a,b}_n = \frac{n!}{\prod_{k=1}^n\left(k \frac{b}{a} -(k-1)\right)},$$ which albeit similar in form is different from the sequence in (\ref{eq:3.25}).

\section{Structural identities}

In this section we will establish commutator relations between the operators and their adjoints.
\subsection{Commutators}
\begin{lemma} \label{Lem:4} For $q >0$ such that  $nq \not= n-1,$ for all $n \in \BN$, we have
\begin{equation} \label{commutatorR_0}
[R_0, R_0^\ast] z^k = \left\{  \begin{array}{cl}
1, & k=0 \\
~& \\
(q-2) (R_0 \mathbb{I}^2 R_0 )z^{k}, & k=1, 2,3,\ldots
\end{array} \right.
\end{equation}
\end{lemma}

\begin{proof}
Looking at the commutator between $R_0$ and $R_0^\ast$ we obtain for the case $k=0$ that
$$[R_0, R_0^\ast] z^0 = R_0 R_0^\ast  1 - R_0^\ast R_0 1 = 1-0 =1,$$
while for the remaining cases we get
\begin{eqnarray*}
[R_0, R_0^\ast]  z^k & = & R_0 R_0^\ast  z^k - R_0^\ast R_0  z^k = \frac{kq -(k-1)}{k+1} R_0 z^{k+1} - R_0^\ast z^{k-1} \\
& = & \frac{kq -(k-1)}{k+1}  z^{k} -  \frac{(k-1)q -(k-2)}{k}  z^{k} = \frac{q-2}{k(k+1)}  z^{k} \\
& = & \frac{q-2}{k(k+1)}  (R_0z^{k+1} ) = (q-2) R_0 \left(\frac{ z^{k+1}}{k(k+1)} \right) \\
& = & (q-2) R_0 \mathbb{I}^2 z^{k-1} = (q-2) (R_0 \mathbb{I}^2 R_0 )z^{k},
\end{eqnarray*}for $k=1,2,3, \ldots$ and $q >0.$
\end{proof}

Thus, we obtain the commutativity between $R_0$ and $R_0^\ast$ in the Hardy space case ($q=2$) for all powers $z^k, k=1,2,\ldots,$ but not for $z^0=1.$

Likewise, for the multiplication operator and its adjoint we have:

\begin{lemma} \label{Lem:4a} For $q >0$ such that  $nq \not= n-1,$ for all $n \in \BN$, we have
\begin{equation} \label{commutatorM_0}
[M_z, M_z^\ast] z^k = \left\{  \begin{array}{cl}
-1, & k=0 \\
~& \\
\frac{q-2}{q} z, & k=1 \\
& \\
(q-2) (\mathbb{I}^\ast)^2  R_0^2 z^{k}, & k=2,3,\ldots
\end{array} \right.
\end{equation}under the convention that $M_z^\ast z^0 = 0.$
\end{lemma}
\begin{proof}For $k=0$ we get
$$[M_z, M_z^\ast] z^0 =  M_z M_z^\ast 1 - M_z^\ast M_z 1 = 0- 1 =-1.$$
For $k=1, 2, \ldots$ we have
\begin{eqnarray*}
[M_z, M_z^\ast]   z^k & = & M_z M_z^\ast z^k - M_z^\ast M_z z^k = M_z \frac{k}{(k-1)q -(k-2)}z^{k-1} - M_z^\ast z^{k+1} \\
 & = &  \frac{k}{(k-1)q -(k-2)}z^{k} - \frac{k+1}{kq -(k-1)}z^{k},
\end{eqnarray*} which leads to
$$[M_z, M_z^\ast] z^1 = \frac{q-2}{q} z^1,$$
and
$$
[M_z, M_z^\ast]   z^k = \frac{q-2}{k(k-1) \big(q - \frac{k-1}{k}\big)\big(q - \frac{k-2}{k-1}\big)  } z^k  =  \frac{q-2}{ (k-1)q - (k-2) }  \frac{z^k}{ kq - (k-1) } $$
$$ = \frac{q-2}{ (k-1)q - (k-2) } \BI^\ast z^{k-1} = (q-2) \BI^\ast \Big(   \frac{z^{k-1}}{ (k-1)q - (k-2) }  \Big) = (q-2) (\BI^\ast)^2 z^{k-2}$$
$$=  (q-2) (\BI^\ast)^2 R_0^2 z^{k}$$
for the remaining values of $k.$
\end{proof}

\subsection{Some $q$-Stirling-like numbers}

One problem one faces is the fact that the structural identities do depend on the powers $z^k$ in which they act, i.e.
$$[M_z, M_z^\ast]   z^k = \lambda(k;q) z^k,\quad k=1, 2, \ldots,$$
whereas
$\lambda(k;q)  = \frac{q-2}{[kq - (k-1)] [(k-1)q - (k-2)]  },$ as seen in Lemma \ref{Lem:4a}.\\

We shall define our $q$-Stirling-like numbers as the coefficients $C_{k}(n, j)$ of the following decomposition
$$(M_z M_z^\ast)^n  = \sum_{j=1}^n C_{k}(n, j) M_z^j (M_z^\ast)^j, \qquad n \in \BN.$$

\begin{lemma} \label{Lem:4.2A}
We have for the $q$-Stirling-like numbers the recursion formula
\begin{eqnarray*}
C_k (1,1) &=& 1; \\
C_{k}(n+1,1)& = & (-\lambda(k;q))^{n} C_k (1,1) \\
& =&  (-\lambda(k;q))^{n}; \\
C_k(n+1, j) & =& C_{k}(n, j-1) - j \lambda(k;q) C_{k}(n, j), \quad j=2,3, \ldots, n;\\
C_k(n+1, n+1) & =&  C_k(n, n) = \cdots = C_k(1,1) =1, 
\end{eqnarray*}
for $n\in \BN.$
\end{lemma}

\begin{proof}
First of all we have the trivial statement $(M_z M_z^\ast)^1=M_z M_z^\ast$ which means $C_k (1,1) = 1$.
Now, assume $(M_z M_z^\ast)^n  = \sum_{j=1}^n C_{k}(n, j) M_z^j (M_z^\ast)^j$. Then
\begin{eqnarray*}
(M_z M_z^\ast)^{n+1} &=& \left( \sum_{j=1}^n C_{k}(n, j) M_z^j (M_z^\ast)^j \right) M_z M_z^\ast \\
&=& \sum_{j=1}^n C_{k}(n, j) M_z^j \left((M_z^\ast)^j M_z \right) M_z^\ast.
\end{eqnarray*}
Due to the commutator relation (\ref{commutatorM_0}) we have $M_z^\ast M_z  = M_z M_z^\ast -\lambda(k;q) = M_z M_z^\ast -\lambda,$ so that
$$(M_z^\ast)^j M_z = (M_z^\ast)^{j-1} ( M_z^\ast M_z) = (M_z^\ast)^{j-1} ( M_z M_z^\ast -\lambda) = (M_z^\ast)^{j-1} M_z M_z^\ast - \lambda (M_z^\ast)^{j-1}$$
$$= (M_z^\ast)^{j-2} ( M_z^\ast M_z) M_z^\ast - \lambda (M_z^\ast)^{j-1} = (M_z^\ast)^{j-2} ( M_z M_z^\ast -\lambda) M_z^\ast - \lambda (M_z^\ast)^{j-1}$$
$$= (M_z^\ast)^{j-2} M_z (M_z^\ast)^2  -2 \lambda (M_z^\ast)^{j-1} = \cdots = M_z (M_z^\ast)^j - j \lambda (M_z^\ast)^{j-1}.$$
Replacing in the initial equation we obtain
\begin{eqnarray*}
(M_z M_z^\ast)^{n+1} &=& \sum_{j=1}^n C_{k}(n, j) M_z^j \left(  M_z (M_z^\ast)^j - j \lambda (M_z^\ast)^{j-1}      \right) M_z^\ast \\
&=& \sum_{j=1}^n C_{k}(n, j) \left( M_z^{j+1} (M_z^\ast)^{j+1} - j \lambda M_z^j ( M_z^\ast)^j \right) \\
&=& \sum_{j=2}^{n+1} C_{k}(n, j-1)  M_z^{j} (M_z^\ast)^{j}  - \lambda \sum_{j=1}^n C_{k}(n, j)  j  M_z^j ( M_z^\ast)^j \\
&=& - \lambda C_{k}(n, 1)  M_z ( M_z^\ast) + \sum_{j=2}^{n} [C_{k}(n, j-1) - j \lambda C_{k}(n, j) ] M_z^j ( M_z^\ast)^j \\
& & \qquad + C_{k}(n, n)  M_z^{n+1} (M_z^\ast)^{n+1}.
\end{eqnarray*}
Comparing the coefficients with the coefficients of $(M_z M_z^\ast)^{n+1}  = \sum_{j=1}^{n+1} C_{k}(n+1, j) M_z^j (M_z^\ast)^j$ we get our recursion formulae.
\end{proof}

For a short example where $k=2$ and $q\not=\frac{1}{2}$ we get
$$\lambda(2, q) = \frac{q-2}{ 2(q-\frac{1}{2})q},$$
while for $k=3$ ($q\not=\frac{1}{2}, \frac{2}{3}$) we have
$$\lambda(3, q) = \frac{q-2}{6 (q-\frac{2}{3})(q-\frac{1}{2})}.$$
Denoting $\lambda:= \lambda(k, q)$ we obtain the following table
\begin{center}
\begin{tabular}{|c||c|c|c|c|c|}
\hline
$C_k(n,j)$ & 1& 2& 3& 4 & ~~5~~ \\
\hline
\hline
1& $1$ &  & & &  \\
\hline
2& $-\lambda$  & $1$ &  &  & \\
\hline
3 & $ \lambda^2$ & $-3\lambda$ & $1$ & & \\
\hline
4 & $ -\lambda^3$ & $7\lambda^2$ & $-6\lambda $ & $1$ & \\
\hline
5 & $ \lambda^4$ & $-15\lambda^3$ & $25\lambda^2$ & $-10\lambda$ & $1$ \\
\hline
\end{tabular}
\end{center}

The coefficients of the power of $\lambda$ in the above table are easily seen to be the signed Stirling numbers of the second kind.

\begin{remark}
For the case of the classic Fock space where $\lambda=-1$ the recursion formulae for the coefficients $C_k(n,j)$ in Lemma~\ref{Lem:4.2A} are the recursion formulae for the Stirling numbers of the second kind~\cite{MR3816055}. In case of the Hardy space we have $\lambda=0$ which corresponds to the fact that the multiplication operator and the backward shift operator are adjoint to each other.
\end{remark}

\section{Calculation of Jordan Chains}

For the Hilbert space $\cH(K_q)$ with reproducing kernel $K_q(z,w) = \qexp{z\overline w}, ~|z\overline w| < \frac{1}{1-q},$ in terms of the $q$-exponential with $q \in (0,1),$ and such that $nq \not=n-1,~n \in \mathbb N,$
\begin{equation} \label{Tsallis_expcalc}
\qexp{x} := \left\{  \begin{array}{cl}
[1+(1-q)x]^{\frac{1}{1-q}}, & 1+(1-q)x \geq 0 \\
 & \\
 0, & \mbox{otherwise}
\end{array}  \right. ,
\end{equation}
we have that it holds $$\lim_{q \to 1} K_q(z,w) = e^{z\overline w},$$
the reproducing kernel for the classical Fock space.

Furthermore, under these conditions we get
\begin{gather}
\qexp{x} := \sum_{k=0}^\infty (1-q)^k \,\Big(\frac{1}{1-q}\Big)_k \frac{(z\ov w)^k}{k !}, \qquad \alpha = 1-q, \nonumber \\
= 1 + z\ov w + \sum_{k=2}^\infty \frac{q (2q-1) \cdots \big[ (k-1)q -(k-2)  \big]}{k !} (z\ov w)^k. \label{Tsallis_expcalc2}
\end{gather}

In the Hilbert space $\cH(K_q)$ we also have
\begin{itemize}
\item the multiplicative operator $M_z z^k = z^{k+1},$ with  adjoint $$M^\ast_z z^k = \frac{k\,z^{k-1}}{(k-1)q -(k-2)},\quad k=0,1, 2, \ldots,$$
\item the backshift operator $R_0 z^0 = 0,$ $R_0 z^k = z^{k-1},$ $k=1, 2, \ldots$, whereas its adjoint is $$R^\ast_0 z^k = [kq -(k-1)] \frac{z^{k+1}}{k+1}, \quad k=0,1, 2, \ldots.$$
\end{itemize}

\begin{lemma}
The adjoint of the multiplicative operator $M^\ast_z$ satisfies
\begin{itemize}
\item[i)] $M^\ast_z \qexp{z} = \qexp{z},$ for all $|z| < \frac{1}{1-q};$
\item[ii)] $\lim_{q\to 1^-} M^\ast_z z^k = \partial_z z^k,$ for all $k =0,1,2,\ldots,$ i.e.  $\lim_{q\to 1^-} M^\ast_z  = \partial_z,$ as operators over the set of polynomials.
\end{itemize}
\end{lemma}

\begin{proof} The first statement is immediate from (\ref{M_z^ast}) with $w=1$. The second statement is obvious.
\end{proof}

We aim to construct Jordan chains based on the existence of an analytic function in $\mathbb D_q : |z| < \frac{1}{|1-q|}$ satisfying to
\begin{equation} \label{GeneralizedEigen}
(M^\ast_z - I) f = \qexp{z}, \qquad z \in \mathbb D_q.
\end{equation}

\begin{theorem} For $q>0$ and such that $nq \not=n-1,$ for all $n\in \mathbb N,$ the above equation (\ref{GeneralizedEigen}) has a solution
\begin{equation} \label{Sol_GeneralizedEigen}
f(z) = f_0+ \sum_{k=1}^\infty  \frac{k+f_0}{k} \frac{(1-q)^{k} \,\Big(\frac{1}{1-q}\Big)_{k}}{(k-1)!} z^k, \qquad z \in \mathbb D_q.
\end{equation}
In particular, for $f_0=0$ we have 
\begin{equation} \label{Sol_GeneralizedEigen1}
f(z) = \sum_{k=1}^\infty   \frac{(1-q)^{k} \,\Big(\frac{1}{1-q}\Big)_{k}}{(k-1)!} z^k, \qquad z \in \mathbb D_q.
\end{equation}
In the latter case we have $f=M_z\partial_z\qexp{z}$ where $\partial_z$ denotes the classic derivative operator.
\end{theorem}

\begin{proof} For a power series $f(z) = \sum_{k=0}^\infty f_k z^k$ and using equation (\ref{GeneralizedEigen}) we get
\begin{gather*}
(M^\ast_z - I) f = \qexp{z} \Leftrightarrow \sum_{k=1}^\infty f_k (M^\ast_z z^k) - \sum_{k=0}^\infty f_k z^k = \qexp{z}\\
\Leftrightarrow \sum_{k=1}^\infty f_k \frac{k\,z^{k-1}}{(k-1)q -(k-2)} - \sum_{k=0}^\infty f_k z^k = \sum_{k=0}^\infty (1-q)^k \,\Big(\frac{1}{1-q}\Big)_k \frac{z^k}{k !}\\
\Leftrightarrow \sum_{k=0}^\infty \left[f_{k+1} \frac{(k+1)}{kq -(k-1)} - f_k \right]z^k = \sum_{k=0}^\infty (1-q)^k \,\Big(\frac{1}{1-q}\Big)_k \frac{z^k}{k !}.
\end{gather*}
This leads to the recursive formula
\begin{equation} \label{Recursive}
f_{k+1} \frac{(k+1)}{kq -(k-1)} - f_k = \frac{(1-q)^k \,\Big(\frac{1}{1-q}\Big)_k}{k !} \quad \Leftrightarrow \quad f_{k+1} = f_k \, \frac{kq -(k-1)}{(k+1)}  + \frac{(1-q)^{k+1} \,\Big(\frac{1}{1-q}\Big)_{k+1}}{(k+1)!},
\end{equation}for $k=0,1,2,\ldots$

Thus, for $k=0$ we get
$f_{1} =  f_0+1,$
while for the remaining terms we obtain
\begin{eqnarray*}
f_{k+1} &=& f_k \, \frac{kq -(k-1)}{(k+1)}  + \frac{(1-q)^{k+1} \,\Big(\frac{1}{1-q}\Big)_{k+1}}{(k+1)!} \\
&=& \left( f_{k-1} \, \frac{(k-1)q -(k-2)}{k}  + \frac{(1-q)^k \,\Big(\frac{1}{1-q}\Big)_k}{k!} \right) \frac{kq -(k-1)}{(k+1)}  + \frac{(1-q)^{k+1} \,\Big(\frac{1}{1-q}\Big)_{k+1}}{(k+1)!} \\
&=& f_{k-1} \, \frac{\big((k-1)q -(k-2) \big)\big( kq -(k-1)\big)}{k(k+1)}    + 2 \frac{(1-q)^{k+1} \,\Big(\frac{1}{1-q}\Big)_{k+1}}{(k+1)!} \\
& \vdots & \\
&=& f_{0} \, \frac{1\cdot q \,(2q-1) \cdots \big((k-1)q -(k-2) \big)\big( kq -(k-1)\big)}{(k+1)!}    + (k+1) \frac{(1-q)^{k+1} \,\Big(\frac{1}{1-q}\Big)_{k+1}}{(k+1)!} \\
&=& f_{0} \,\frac{(1-q)^{k+1} \,\Big(\frac{1}{1-q}\Big)_{k+1}}{(k+1)!}  + \frac{(1-q)^{k+1} \,\Big(\frac{1}{1-q}\Big)_{k+1}}{k!}.
\end{eqnarray*}
Hence, the solution to equation (\ref{GeneralizedEigen}) has the form
$$f (z) = \sum_{k=0}^\infty f_k z^k = f_0+\sum_{k=1}^\infty  \left[  f_{0} \,\frac{(1-q)^{k} \,\Big(\frac{1}{1-q}\Big)_{k}}{k!}  + \frac{(1-q)^{k} \,\Big(\frac{1}{1-q}\Big)_{k}}{(k-1)!} \right] z^k,$$
and easy calculations show this series to converge for $z \in \mathbb D_q : |z| < \frac{1}{|q-1|}.$

Choosing $f_0=0$ we have 
$$
f (z) = \sum_{k=1}^\infty f_k z^k = \sum_{k=1}^\infty  \left[  \frac{(1-q)^{k} \,\Big(\frac{1}{1-q}\Big)_{k}}{(k-1)!} \right] z^k.
$$
Now observing that $M_z\partial_z z^k=kz^k$ we get $f=M_z\partial_z \qexp{z}$.
\end{proof}

Let us remark that in the case of $q\rightarrow 1$ we have $f(z)=ze^z$ while in the case of $q=2$ we obtain $f(z)=\frac{-z}{(1-z)^2}$.

Using $f(\lambda z)$ instead of $f(z)$ in the above proof leads to the following corollary.

\begin{corollary}
For $q>0$ and such that $nq \not=n-1,$ for all $n\in \mathbb N,$ the equation 
$$
(M^\ast_z - \lambda I) f_\lambda(z) = \qexp{\lambda z}, \qquad z \in \mathbb D_q.
$$
has a solution
\begin{equation} \label{Sol_GeneralizedEigen1a}
f_\lambda(z) = f_0 +\sum_{k=1}^\infty  \frac{k+f_0}{k} \frac{(1-q)^{k} \,\Big(\frac{1}{1-q}\Big)_{k}}{(k-1)!} \lambda^k z^k, \qquad z \in \mathbb D_q.
\end{equation}
In particular, for $f_0=0$ we have 
\begin{equation} \label{Sol_GeneralizedEigen1b}
f_\lambda (z) = \sum_{k=1}^\infty   \frac{(1-q)^{k} \,\Big(\frac{1}{1-q}\Big)_{k}}{(k-1)!} \lambda^k z^k, \qquad z \in \mathbb D_q.
\end{equation}

\end{corollary}

The previous corollary allows us to construct Jordan chains in the setting of Tsallis.

\section{Rational functions in the sense of Tsallis}
\setcounter{equation}{0}

Let us consider matrix-valued rational functions in this context. To this end we recall the definition of rational functions related to a generalized derivative operator $\partial_\varphi$.

\begin{definition}
$f$ is rational if the linear span of the functions $f,\partial_\varphi f,\partial_\varphi^2f,\ldots$ is finite dimensional.
 \end{definition}

In case of the backward shift operator $\partial_\varphi=R_0$ we have the classic definition of a matrix-valued rational function.\\

In the case of the classic derivative operator $\partial_\varphi=\partial_z$ we get rational functions as matrix-valued functions of the form
   \[
     f=Ce^{zA}B,
   \]
  for details, see~\cite{daf1}.
As we will see in the case of Tsallis $q$-calculus with $\partial_\varphi=M_z^*$ we have rational functions to be of the form
     \[
     f=Ce_q^x(zA)B.
   \]
   In the sequel we are going to prove the equivalence of different characterizations of matrix-valued rational functions in the sense of Tsallis. To this end we introduce the Tsallis version of a $q$-Borel transform.

 \begin{definition}
Given a function
$$
F(z)=\sum_{k=0}^\infty F_k z^k
$$
analytic in the disk $|z|<r$ its $q$-Tsallis Borel transform $F \to B_qF$ is defined as the function
 $$
B_qF(z)=  \sum_{k=0}^\infty \frac{(1-q)^k \,\Big(\frac{1}{1-q}\Big)_k}{k!} F_k z^k.
 $$

 \end{definition}

 \begin{lemma} If the function $
F(z)=\sum_{k=0}^\infty F_k z^k
$ is analytic in the disk $|z|<r$ then its  $q$-Tsallis Borel transform $B_qF$ is analytic in the disk $|z|<\frac{r}{|1-q|}$.
 \end{lemma}

 The proof of this lemma is immediate keeping in mind that using (\ref{RKHS})  the series  $$
 \sum_{k=0}^\infty \frac{(1-q)^k \,\Big(\frac{1}{1-q}\Big)_k}{k!}F_k z^k
 $$
 we have the radius of convergence $\frac{r}{|1-q|}$ if the series $\sum_{k=0}^\infty F_k z^k$ has radius of convergence $r$.

Since the $q$-Tsallis Borel transform defines a one-to-one mapping from this lemma we get the following characterization of a $q$-rational function in the sense of Tsallis.

\begin{theorem}
A function $F$ is a $q$-rational matrix-valued function in the sense of Tsallis if and only if it is the  $q$-Tallis Borel transform of a classical rational matrix-valued function.
\end{theorem}

 This leads to the following characterizations of matrix-valued rational functions in the sense of Tsallis.

 \begin{theorem}
  For a power series $F(z)=\sum_{k=0}^\infty F_k z^k$  converging in a neighborhood of the origin, the following characterizations of $q$- rationality
  are equivalent:
\begin{enumerate}
\item {\bf Realization.} $F$ admits a realization, i.e. can be written as
\[
F(z)=C\left(e_q^x(zA)\right)B,
\]
with
where $N\in\mathbb N_0$ and $(C, A,B)\in \mathbb C^{n\times N} \times \mathbb C^{N\times N}\times \mathbb C^{N\times m}$. If $N=0$, $R(z)=CB$.\\

\item {\bf Hankel operator.} The block Hankel operator
\[
H=\begin{pmatrix}{F_0}  &{F_1} &{F_2}\frac{2!}{q} &\udots &\udots\\
 & & & & \\
   {F_1} & {F_2}\frac{2!}{q} & {F_3}\frac{3!}{q(2q-1)}&\udots &\udots\\
    & & & & \\
    {F_2}\frac{2!}{q}& {F_3}\frac{3!}{q(2q-1)}&\udots &\udots &\udots\\
            & & & &  \\
  \udots  &\udots & \udots & {F_{2k}}\frac{(2k)!}{q(2q-1)\cdots ((2k-1)q -(2k-2))} &\udots \\
              & & & &  \\
  \udots  &\udots & \udots &\udots &\udots
  \end{pmatrix}
  \]
  has finite rank.\\

\item {\bf Taylor series.}
The coefficients $F_k$ can be written as
\[
F_k=\frac{ (1-q)^k \,\Big(\frac{1}{1-q}\Big)_k }{k!} CA^{k}B,\quad k=0,1,2,\ldots
\]
for some matrices $(C, A,B)\in \mathbb C^{n\times N} \times \mathbb C^{N\times N}\times \mathbb C^{N\times m}$. 
If $N=0$, all the $F_k=0$, $k=1,2,\ldots$.\\

\item {\bf Characterization by the adjoint of $M_z$} With $M_z^*$ as in (\ref{M^ast}), the linear span of the columns of the functions $(M_z^*)^kF$, $k=0,1,2\ldots$ is finite dimensional.\\

\end{enumerate}
\end{theorem}
For the proof we remark that the realization $(1)$ corresponds to the $q$-Tsallis Borel transform of a classical rational matrix-valued function. Furthermore, $(3)$ is $q$-Tsallis Borel transform of the corresponding Taylor series of the corresponding classical rational matrix-valued function. For the last item we remark that we have $M_z^*T=R_0B_qT$ where $R_0$ is the classic backward-shift operator.

Furthermore, we can state the following theorem which represents an adaptation of Theorem~9.10 in~\cite{daf1} (page 423). 

\begin{theorem}
The following statements are equivalent:
\begin{enumerate}
\item The columns of $M_z^*F$ belong to a finite-dimensional $M_z^*$-invariant space.
\item $F(z)$ is of the form
$$
F(z)=D+I_z (Ce_q^x(zA)B)
$$
where $I_z$ is the q-Jackson integration operator.
\end{enumerate}

\end{theorem}

Let us remark that the above approach is a special case of a general setting connected to derivative operators of Gelfond-Leontiev type. Let us recall its definition, see, e.g.,~\cite{Gelfond51, MR1265940}.

\begin{definition}
Let
\begin{equation} \label{Eq:EntireFunction}
\varphi(z) =\sum_{k=0}^\infty \varphi_k z^k,
\end{equation}
be an entire function with order $\rho >0$ and  degree $\sigma > 0,$ that is, such that $\lim_{k \rightarrow \infty} k^{\frac{1}{\rho}} \sqrt[k]{|\varphi_k|} =\left(\sigma e \rho\right)^{\frac{1}{\rho}}.$ We define the Gelfond-Leontiev (G-L) operator of generalized differentiation with respect to $\varphi,$ denoted as $\partial_\varphi,$
as the operator acting on an analytic function $f(z) =\sum_{k=0}^\infty a_k z^k$ as
\begin{equation} \label{Eq:p1+p2}
f(z) =\sum_{k=0}^\infty a_k z^k \quad \mapsto \quad \partial_\varphi f(z) =\sum_{k=1}^\infty a_k \frac{\varphi_{k-1}}{\varphi_k} ~z^{k-1}. \end{equation}
\end{definition}

In this case the corresponding Borel transform has the form $B_\varphi F$ 
 $$
\sum_{k=0}^\infty F_k z^k\to B_\varphi F (z) = \sum_{k=0}^\infty \varphi_k F_k z^k
 $$
and with the same argument as above we have that following theorem.

\begin{theorem}
  For a power series $F(z)=\sum_{k=0}^\infty F_k z^k$  converging in a neighborhood of the origin, the following characterizations of rationality are equivalent:
\begin{enumerate}
\item {\bf Realization.} $F$ admits a realization, i.e. can be written as
\[
F(z)=C\left(\varphi(zA)\right)B,
\]
with
where $N\in\mathbb N_0$ and $(C, A,B)\in \mathbb C^{n\times N} \times \mathbb C^{N\times N}\times \mathbb C^{N\times m}$. If $N=0$, $R(z)=CB$.\\

\item {\bf Hankel operator.} The block Hankel operator
\[
H=\begin{pmatrix}{F_0}  &{F_1} &{F_2}\varphi_2 &\udots &\udots\\
 & & & & \\
   {F_1} & {F_2}\varphi_2 & {F_3}\varphi_3&\udots &\udots\\
    & & & & \\
    {F_2}\varphi_2& {F_3}\varphi_3&\udots &\udots &\udots\\
            & & & &  \\
  \udots  &\udots & \udots & {F_{2k}}\varphi_k&\udots \\
              & & & &  \\
  \udots  &\udots & \udots &\udots &\udots
  \end{pmatrix}
  \]
  has finite rank.\\

\item {\bf Taylor series.}
The coefficients $F_k$ can be written as
\[
F_k=\varphi_k CA^{k}B,\quad k=0,1,2,\ldots
\]
for some matrices $(C, A,B)\in \mathbb C^{n\times N} \times \mathbb C^{N\times N}\times \mathbb C^{N\times m}$. 
If $N=0$, all the $F_k=0$, $k=1,2,\ldots$.\\

\item {\bf Characterization by the adjoint of $M_z$} The linear span of the columns of the functions $\partial_\varphi^k F$, $k=0,1,2\ldots$ is finite dimensional.\\

\end{enumerate}
\end{theorem}

Furthermore, in the same way as above we can state the following theorem:

\begin{theorem}
The following statements are equivalent:
\begin{enumerate}
\item The columns of $M_z^*F$ belong to a finite-dimensional $M_z^*$-invariant space.
\item $F(z)$ is of the form
$$
F(z)=D+I_z (C\varphi(tA)B)
$$
where $I_z$ is the integration operator of Gelfond-Leontiev type with respect to $\varphi$.
\end{enumerate}
\end{theorem}

\section{Aknowledgements}
Daniel Alpay thanks the Foster G. and Mary McGaw Professorship in Mathematical Sciences, which supported this research. The second and third author were supported by CIDMA, through the Portuguese FCT (UIDP/04106/2025 and UIDB/04106/2025). 


     \bibliographystyle{plain}
\def\cprime{$'$} \def\cprime{$'$} \def\cprime{$'$}
  \def\lfhook#1{\setbox0=\hbox{#1}{\ooalign{\hidewidth
  \lower1.5ex\hbox{'}\hidewidth\crcr\unhbox0}}} \def\cprime{$'$}
  \def\cprime{$'$} \def\cprime{$'$} \def\cprime{$'$} \def\cprime{$'$}
  \def\cprime{$'$}

   \end{document}